\documentclass[11pt]{amsart}

\usepackage{epsfig}
\usepackage{amsthm}
\usepackage{amssymb}
\usepackage{amsmath}
\usepackage{epic}
\usepackage{eepic}
\usepackage{graphicx}

\newcommand     {\comment}[1]   {}
\newcommand{\mute}[2] {}
\newcommand     {\printname}[1] {}
\newtheorem {Theorem}   {Theorem}
\theoremstyle{definition}

\theoremstyle{remark}

\newtheorem {Corollary}{Corollary}

\def    \bfC    {{\mathbf C}}

\def    \bfone  {{\mathbf 1}}

\def    \C      {\mathbb{C}}
\def    \R      {\mathbb{R}}

\def    \K      {\mathcal{K}}

\def    \codim    {{\operatorname{codim}}}

\newcommand{\dsum}{\displaystyle\sum}  % Sum in displaystyle
\begin{document}

\bibliographystyle{plain}

\title[Weighted Brianchon-Gram decomposition]
{Weighted Brianchon-Gram decomposition}

\author[J.\ Agapito]{Jos\'e Agapito}
\address{DEPARTAMENTO DE MATEM\'ATICA, INSTITUTO SUPERIOR T\'ECNICO, AV. ROVISCO PAIS,
1049-001 LISBON, PORTUGAL, FAX: (351) 21 841 7035}
\email{agapito@math.ist.utl.pt}

\thanks{2000 \emph{Mathematics Subject Classification.}
Primary 52B}
\thanks{Partially supported by NSF (United States) through grants DMS 99/71914 and DMS 04/05670 
and by FCT (Portugal) through program POCTI/FEDER and grant POCTI/SFRH/BPD/20002/2004.}

\begin{abstract}
We give in this note a weighted version of Brianchon-Gram's decomposition for a simple
polytope. We can derive from this decomposition the weighted polar formula of \cite{A} and a weighted
version of Brion's theorem \cite{B88}, in a manner similar to \cite{H}, where the unweighted case is
worked out. This weighted version of Brianchon-Gram is a direct consequence of the ordinary Brianchon-Gram
formula.
\end{abstract}

\maketitle

%\tableofcontents

% -------------------------------------------------------------------------
\section{Introduction}
% -------------------------------------------------------------------------
\label{se:intro} Let $A\subset\R^n$ be a closed convex subset. The characteristic function
$\bfone_A$ of $A$ is the function $\bfone_A\colon\R^n\to\C$ given by
\begin{equation*}
\bfone_A(x)=\left\{\begin{array}{cc}
                        1 & \mbox{if }x\in A \\
                        0 & \mbox{if }x\notin A \\
                      \end{array}\right. .
\end{equation*}

Let $\K(\R^n)$ be the complex vector space spanned by the functions $\bfone_A$. Thus, a
function $f\in\K(\R^n)$ is a linear combination
$$f=\sum_{i=1}^m \alpha_i\bfone_{A_i},$$
where the $A_i$ are closed convex sets in $\R^n$ and the $\alpha_i$ are complex numbers.\par

Among the elements of $\K(\R^n)$ there are three well known decomposition formulas: the
Brianchon-Gram decomposition \cite{Br}, \cite{G} (see also \cite{B} and \cite{S}), which
determines the characteristic function of any polytope as a signed sum of characteristic
functions of cones associated to its faces, the polar decomposition of a simple polytope $P$
(\cite{L}, \cite{V}), which uses the notion of polarization and Brianchon-Gram's formula 
in order to write the characteristic function of $P$ in terms of the characteristic functions 
of cones based on the vertices of $P$ only, and the Brion decomposition of a polytope \cite{B88}, 
which is also a direct consequence of Brianchon-Gram's formula.\par

\begin{figure}[h]
  \centering
  \includegraphics[scale=.85]{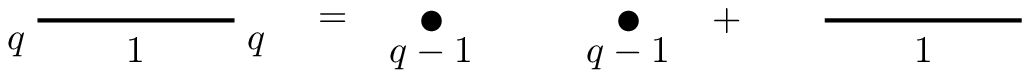}\\
  \caption{Decomposition of an interval.}\label{fi:interval-cell}
\end{figure}

We can put \emph{weights} to the faces of $P$ in a meaningful way and get a new element of
$\K(\R^n)$. For instance, let $q$ be any complex number and let $[a,b]$ be any interval. We can
write (see Figure \ref{fi:interval-cell})

\noindent That is, as an element of $\K(\R)$, the function
$(q-1)\bfone_{\{a\}}+(q-1)\bfone_{\{b\}}+ \bfone_{[a,b]}$ defines a
weighted characteristic function over $[a,b]$. We denote it by
$\bfone^q_{[a,b]}$. In general, we can assign arbitrary complex
values to the facets of a polytope $P$ and construct a new element
$\bfone^w_P$ of $\K(\R^n)$. (See \eqref{q-Delta}.)
\begin{figure}[h]
  \centering
  \includegraphics[scale=.85]{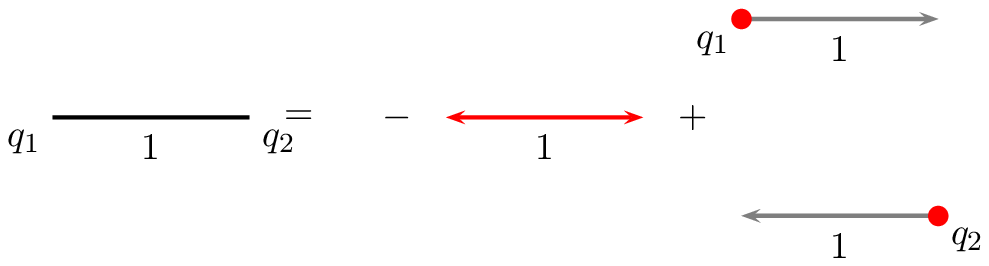}\\
  \caption{Weighted Brianchon-Gram decomposition for an interval.}\label{fi:q-interval}
\end{figure}

When we assign the same value $q$ to the facets of $P$, we also denote $\bfone^w_P$ by
$\bfone^q_P$. The weighted polar decomposition formula of \cite{A} expresses $\bfone^q_{P}$ as
an alternating sum of weighted characteristic functions of cones based on the vertices of $P$.
For example, in the case of $\bfone^q_{[a,b]}$, we have
$\bfone^q_{[a,b]}=\bfone^q_{[a,\infty)}-\bfone^{1-q}_{(-\infty,b]}$. The purpose of this note
is to show that $\bfone^w_P$ (defined in \eqref{q-Delta}) satisfies a weighted version of
Brianchon-Gram's formula, from which it readily follows the weighted polar formula of \cite{A}
and a weighted version of Brion's theorem \cite{B88}. The relationship among these formulas has
already been pointed out by Haase \cite{H} in the unweighted case. Our main result is the
weighted version of Brianchon-Gram's formula as stated in \eqref{BrianchonGram}. (See Figures
\ref{fi:q-interval}, \ref{fi:q-triangle} and \ref{fi:w-pyramid} for illustrations of this
formula.)

% -------------------------------------------------------------------------
\section{The weighted formula}
% -------------------------------------------------------------------------
\label{se:weighted-formula}

Let $P$ be a $d$-dimensional polyhedron in $\R^d$ (for standard definitions on polyhedra we
refer to \cite{B}). We can write it as the intersection of a finite number of half-spaces
\begin{equation}\label{eq:polyhedron}
P = H_1\cap\ldots\cap H_N ,
\end{equation}
where $H_i=\{x\,|\,\langle u_i,x\rangle + \mu_i\ge 0\}$, with $\mu_i\in\R$ and $u_i\in(\R^d)^*$
for $1\le i\le N$. Note that $\R^d = \{x\,|\,\langle 0,x\rangle + 0 \ge 0\}$, hence $\R^d$ is
trivially a polyhedron. It follows that $P$ is a closed convex set. We assume that $P$ is
obtained with the smallest possible $N$. The facets of $P$ are $\sigma_i=P\cap\partial H_i$ for
$i=1,\ldots,N$. If the intersection \eqref{eq:polyhedron} is bounded then $P$ is a polytope. We
say that a $d$-dimensional polyhedron $P$ is simple if every vertex of $P$ belongs (when it
exists) to exactly $d$ facets of $P$. In the case of a polyhedron without vertices, we assume
that this condition is trivially satisfied.

Let $P$ be any $d$-dimensional polyhedron in $\R^d$. For each $i=1,\ldots,N$, we assign
arbitrary complex numbers $q_i$ to the facets $\sigma_i$ of $P$. Each non-trivial face $F$ of
$P$ ($F\neq\phi,P$) can be uniquely described as an intersection of facets
\begin{equation}\label{eq:face}
F = \bigcap_{i\in I_F}\sigma_i,
\end{equation}
where $I_F$ denotes the set of all facets of $P$ containing $F$. When $P$ is simple, the number
of elements in $I_F$ is equal to the codimension of $F$.

To each non-trivial face $F$ we assign the value $\prod_{i\in I_F}q_i$. When $F=P$, we give it
the value 1. This amounts to defining the weighted function $w\colon P\to\C$ by
$w(x)=\prod_{i\in I_F}q_i$, where $F$ is the face of $P$ of smallest dimension containing $x$.
If $x$ is in the interior of $P$, we set $w(x)=1$. We extend this definition to all $\R^d$ and
get the weighted characteristic function
\begin{equation}\label{q-Delta}
\bfone^w_{P}(x)=\left\{\begin{array}{cc}
                        w(x) & \mbox{if }x\in P \\
                        0 & \mbox{if }x\notin P \\
                      \end{array}\right. .
\end{equation}

\noindent Now, let $F$ be any face of $P$. The tangent cone to $P$ at $F$ is
\begin{equation*}
\bfC_F=\{y+r(x-y)\,|\, r\ge0,\, y\in F,\, x\in P\}.
\end{equation*}
It follows that $\bfC_F$ is also a polyhedron. For example, when $F=P$, we have $\bfC_F=\R^d$.
\begin{figure}[h]
  \centering
  \includegraphics[scale=.85]{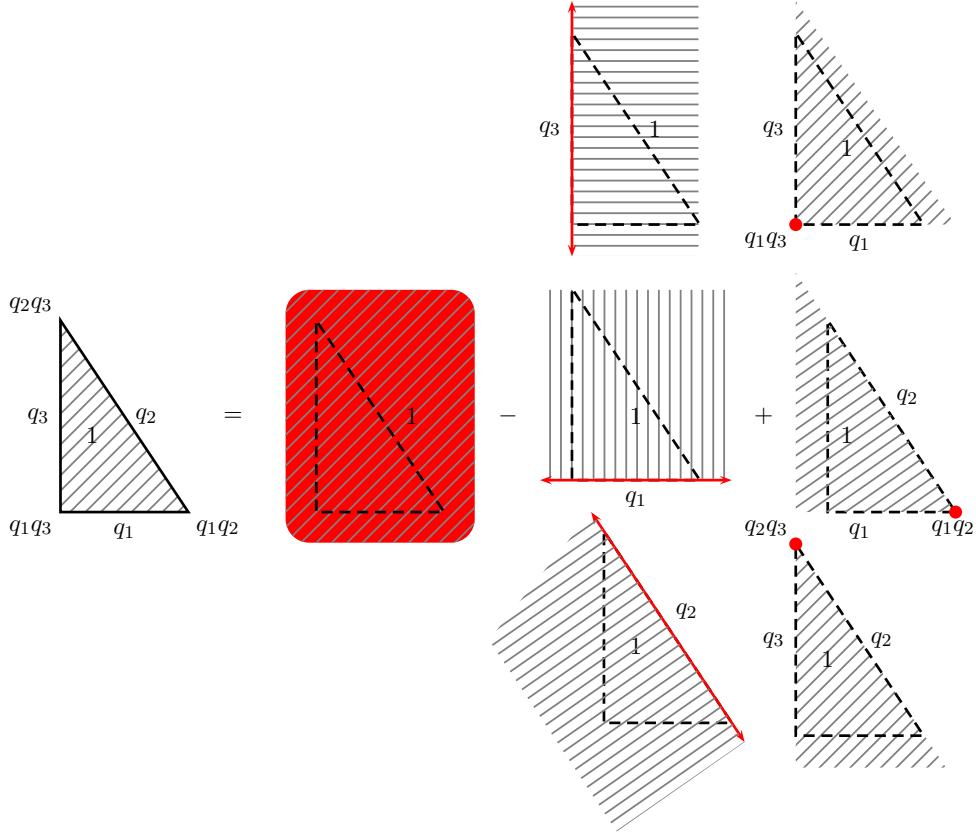}\\
  \caption{Weighted Brianchon-Gram decomposition for a triangle.}\label{fi:q-triangle}
\end{figure}
\begin{Theorem}[Weighted Brianchon-Gram]\label{th:BrianchonGram}
Let $P$ be a simple polytope of dimension $d$ in $\R^d$. We have
\begin{equation} \label{BrianchonGram}
\bfone^w_{P} = \sum_{F\preceq P}(-1)^{\dim F}\bfone^w_{\bfC_F},
\end{equation}
where the sum is over all faces $F$ of $P$.
\end{Theorem}

When $q_1=\ldots=q_N=1$, we have the ordinary Brianchon-Gram formula. We illustrate this
theorem for a triangle in Figure \ref{fi:q-triangle}. (See also Figure
\ref{fi:q-interval}.)\par

\begin{proof}
Let $\Sigma$ be the set of facets of $P$ and let $F$ be a proper face of $P$. Let
$I_F=\{\sigma\in\Sigma\,|\,F\subset\sigma\}$. Then $F=\bigcap_{i\in I_F}\sigma_i$. Since $P$ is
simple, the cardinality of $I_F$ is equal to the codimension of $F$ and for all possible
non-empty subsets $J$ of $I_F$, the face $\bigcap_{i\in J}\sigma_i$ of $P$ contains $F$. A
straightforward computation shows that
\begin{equation}\label{eq:computation}
\prod_{i\in I_F} q_i = 1 + \sum_{\phi\neq J\subset I_F}\,\,\prod_{i\in J} (q_i-1),
\end{equation}
where the $q_i$ are arbitrary complex numbers assigned to the facets $\sigma_i$ of $P$. We
decompose $P$ into all its faces $F$ (including $F=P$). By \eqref{q-Delta} and
\eqref{eq:computation}, we have
\begin{equation}\label{eq:polytope-faces}
\bfone^w_P = \bfone_P + \sum_{F\neq P} \prod_{i\in I_F} (q_i - 1) \bfone_F\, .
\end{equation}
We can apply the ordinary Brianchon-Gram formula to each $F$ and get
\begin{equation}\label{eq:usingBG}
\bfone^w_P = \sum_{F\preceq P}(-1)^{\dim F}\bfone_{\bfC_F} + \sum_{F\neq P} \prod_{i\in I_F} (q_i
-1)\sum_{G\preceq F} (-1)^{\dim G} \bfone_{\bfC_G}\, .
\end{equation}
Note that \eqref{eq:polytope-faces} also holds for the tangent cones $C_F$ of $P$ at $F$; that
is
\begin{equation}\label{eq:cone-faces}
\bfone^w_{\bfC_F} = \bfone_{\bfC_F} + \sum_{H\neq \bfC_F} \prod_{i\in I_H} (q_i - 1) \bfone_H\,
.
\end{equation}
Since the faces $H$ of the tangent cone $\bfC_F$ are in turn tangent cones associated to the
faces $G$ of $F$, we can regroup the characteristic functions in \eqref{eq:usingBG} according
to \eqref{eq:cone-faces} and obtain
\begin{equation*}
\bfone^w_{P} = \sum_{F\preceq P}(-1)^{\dim F}\bfone^w_{\bfC_F}.
\end{equation*}
\end{proof}

Theorem \ref{th:BrianchonGram} can be extended to non-simple
polytopes where the only faces $F$ which are \emph{non-generic}
(i.e. $F$ has dimension $f$ but $\vert I_F\vert\neq d-f$) are
vertices.

\begin{figure}[h]
  \centering
  \includegraphics[scale=.90]{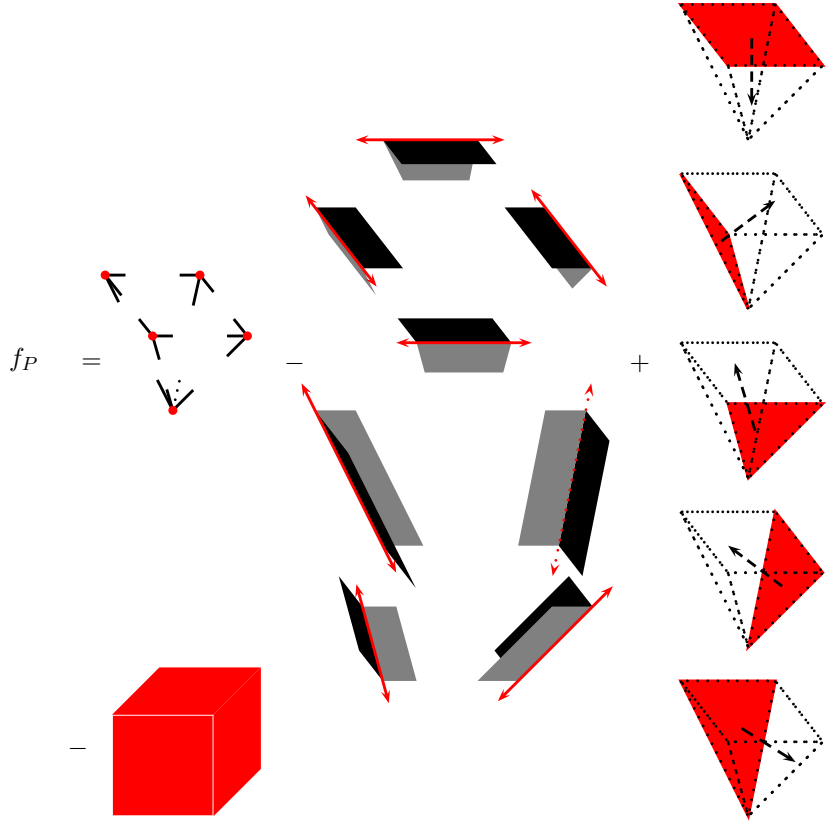}\\
  \caption{Weighted Brianchon-Gram decomposition for a pyramid.}\label{fi:w-pyramid}
\end{figure}

\begin{figure}[h]
  \centering
  \includegraphics[scale=.90]{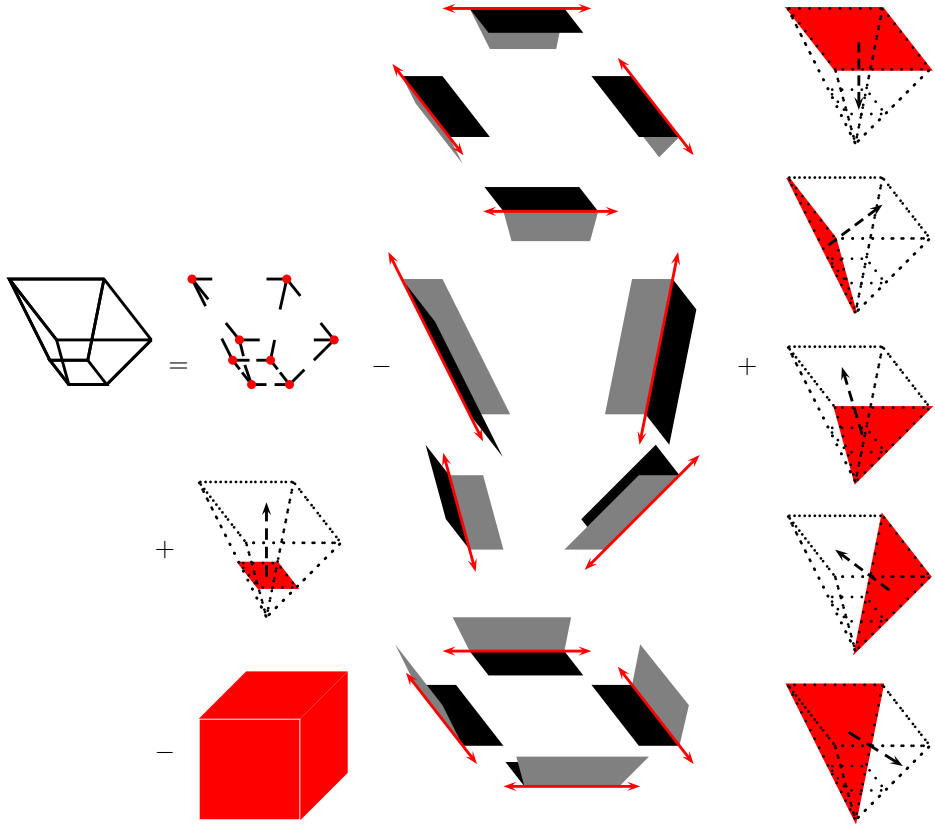}\\
  \caption{Weighted Brianchon-Gram for a truncated pyramid.}\label{fi:w-truncpyramid}
\end{figure}

\begin{figure}[h]
  \centering
  \includegraphics[scale=.85]{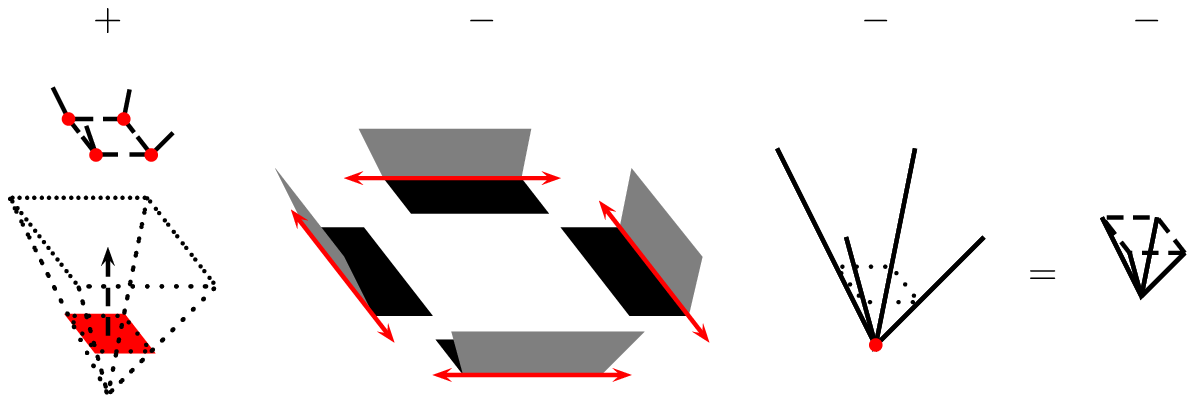}\\
  \caption{Key difference $\bfone^w_{P_s} - f_P$ for a pyramid.}\label{fi:keydifference}
\end{figure}

\begin{figure}[h]
  \centering
  \includegraphics[scale=1]{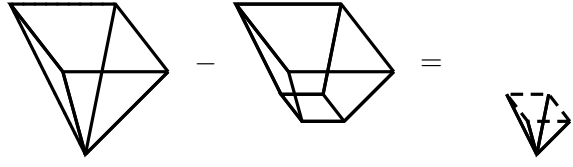}\\
  \caption{Difference of the pyramid and the truncated pyramid.}\label{fi:w-pyr-trunc}
\end{figure}

\noindent Indeed, in dimension three, only vertices may be non-generic for an arbitrary
non-simple polytope. In higher dimensions though, there may be other faces which are the
intersection of \emph{too many} facets. We illustrate this special extension for a
pyramid in Figures \ref{fi:w-pyramid},\ref{fi:w-truncpyramid},\ref{fi:keydifference}
and \ref{fi:w-pyr-trunc}.

Let $P$ be any non-simple polytope of dimension $d$ in $\R^d$ whose
non-generic faces are only vertices. Let $V_{ns}(P)$ be the set of
non-simple vertices of $P$. If $v\in V_{ns}(P)$, then $\vert
I_{v}\vert
> d$. (Recall that in the simple case $\vert I_v\vert = d$ for all vertices of $P$.)

\noindent We chop off all the non-simple vertices $v$ of $P$ by taking
hyperplanes $\sigma_{v}$ very close to these vertices. We orient the
$\sigma_{v}$ away from $v$ and denote the corresponding half-space
associated to $\sigma_{v}$ by $H_{v}$. We assign the constant value
1 to the hyperplanes $\sigma_{v}$. We obtain a simple polytope $P_s$
for which formula \eqref{BrianchonGram} holds; that is
\begin{equation}\label{choppolytope}
\bfone^w_{P_s} = \sum_{F_s\preceq P_s}(-1)^{\dim F_s}\bfone^w_{\bfC_{F_s}},
\end{equation}
where the sum is over all the faces $F_s$ of $P_s$. On the other
hand, we set
\begin{equation}\label{nspolytope}
f_P = \sum_{F\preceq P}(-1)^{\dim F}\bfone^w_{\bfC_{F}},
\end{equation}
where the sum is now over all the faces $F$ of the non-simple polytope
$P$.

\noindent We can clearly see that $\bfone^w_{\bfC_{F}} =
\bfone^w_{\bfC_{F_s}}$ for all the faces of $P$ and $P_s$ which are not
the non-simple vertices $v\in V_{ns(P)}$ and are not the faces of
$P_s$ contained in the various $\sigma_{v}$. Thus, we have
\begin{equation}\label{keydifference}
\bfone^w_{P_s} - f_P = \dsum_{v\in
V_{ns}(P)}\left(\dsum_{F_s\subset\sigma_v}(-1)^{\dim
F_s}\bfone^w_{\bfC_{F_s}}-\bfone^w_{\bfC_v}\right).
\end{equation}

\noindent Notice that all the cross sections of $\bfC_{v}$ contained in $H_{v}$
and parallel to $\sigma_{v}$, are simple polytopes (of the same
type) in $\R^{d-1}$. This is a consequence of the imposed condition
on $P$, that its only non-generic faces be vertices. We can apply
the weighted Brianchon-Gram formula to these cross sections and get
weighted characteristic functions over its faces, which are equal in
absolute value but with opposite signs, to the weighted
characteristic functions over the faces of the corresponding cross
sections of the cones $\bfC_{F_s}$, where $F_s\subset\sigma_v$, for
all $v\in V_{ns}(P)$. Then, we can write \eqref{keydifference} as
\begin{equation}
\bfone^w_{P_s} - f_P = - \sum_{v\in V_{ns}(P)}
\bfone^w_{\bfC_{v}\setminus H_{v}}.
\end{equation}

\noindent On the other hand, it can be easily checked that
\begin{equation}
\bfone^w_{P_s} - \bfone^w_{P} = - \sum_{v\in V_{ns}(P)}
\bfone^w_{\bfC_{v}\setminus H_{v}}.
\end{equation}

\noindent Therefore, we conclude that
\begin{equation}
\bfone^w_{P} = \sum_{F\preceq P}(-1)^{\dim F}\bfone^w_{\bfC_{F}}.
\end{equation}

As an immediate consequence of Theorem \ref{th:BrianchonGram} we obtain a weighted version of
Brion's theorem\footnote{Brion already proved a more general weighted version of his formula in
\cite{B88} (see page 82 of \cite{B92})} \cite{B88}. (See also \cite{B92}.)

\begin{Corollary} For any simple polytope $P$, we have
\begin{equation*} \label{Brion}
\bfone^w_{P} = g + \sum_{v}\bfone^w_{\bfC_{v}},
\end{equation*}
where $g$ is a linear combination of characteristic functions of
cones with straight lines and the sum is over all vertices $v$ of
$P$.
\end{Corollary}
This readily follows from grouping together all tangent cones in \eqref{BrianchonGram} that
contain straight lines, just as in \cite{H}.

When $P$ is a simple polyhedron and the same value $q\in\C$ is assigned to all its facets, the
weighted characteristic function \eqref{q-Delta} gets the form
\begin{equation}
\bfone^w_{P}(x) = \bfone^q_P(x) = \left\{\begin{array}{cc}
                        q^{\codim(F)} & \mbox{if }x\in P \\
                        0 & \mbox{if }x\notin P \\
                      \end{array}\right. ,
\end{equation}
where $F$ is the face of $P$ of smallest dimension containing $x$. Theorem
\ref{th:BrianchonGram} implies the weighted polytope decomposition of \cite{A}. To show this,
let $\xi$ be a polarizing vector\footnote{a polarizing vector is a
generic element of $(\R^n)^*$ which is nonconstant on each edge of $P$. (See \cite{L}, \cite{V}
and compare with \cite{A}.)} in $(\R^d)^*$ and let $\Delta$ be a simple polytope in $\R^d$.
(We now follow the definitions and notation of \cite{A}.) We put together the faces of $P$
according to where they achieve their minimum in the $\xi$-direction. We obtain
\begin{equation*}\label{brianchongram-polar}
(-1)^{\#v} \bfone^{w_v}_{\bfC^\sharp_v} = \sum_{
\substack{v \preceq F \preceq P \\[.5mm] \xi(v) \le \xi(F)}}
(-1)^{\dim F} \bfone^q_{\bfC_F},
\end{equation*}
\noindent where $\#v$ denotes the number of edges of $\bfC_v$
\emph{flipped} according to $\xi$ and where
$\bfone^{w_v}_{\bfC^\sharp_v}$ is the weighted characteristic
function of the $\xi$-polarized tangent cone $\bfC^\sharp_v$
defined in \cite{A}. (To agree with the weighted formulas in
\cite{A} we use the substitution $q=1/(1+y)$.) Then
\begin{equation*}
\sum_{v}(-1)^{\#v} \bfone^{w_v}_{\bfC^\sharp_v} = \sum_{v}\sum_{
\substack{v \preceq F \preceq P \\[.5mm] \xi(v) \le \xi(F)}}
(-1)^{\dim F} \bfone^q_{\bfC_F}=\sum_{F}(-1)^{\dim F}\bfone^q_{\bfC_F}=\bfone^q_{P}.
\end{equation*}

Thus, the weighted polar decomposition of \cite{A} is a direct consequence of Theorem
\ref{th:BrianchonGram}, which can be proved exactly in the same fashion as in \cite{L}.

\end{document}